\documentclass[12pt]{amsart}
\usepackage{amssymb, amsmath, url, graphicx}

\theoremstyle{definition}

\newcommand{\0}{\leqno}

\title{A nilpotency criterion for finite groups}

\author{Marius T\u arn\u auceanu}
\address{Marius T\u arn\u auceanu \\
Faculty of  Mathematics \\
''Al.I. Cuza'' University \\
Ia\c si \\
Romania}
\email{tarnauc@uaic.ro}

\date{2016/2017}

\begin{document}

\begin{abstract}
    Let $G$ be a finite group. In this short note, we give a criterion of nilpotency of $G$ based on the existence of elements of certain order in each section of $G$.
\end{abstract}


\maketitle

\section{Introduction}

The problem of detecting structural properties of finite groups by looking at element orders has been considered in many recent papers (see e.g. \cite{1} and \cite{3}-\cite{6}).
In the current note, we identify a new property detecting nilpotency of a finite group $G$ that uses the function
$$\varphi(G)=|\{a\in G \mid o(a)=\exp(G)\}|$$introduced and studied in \cite{9}. The proof that we present is founded on the structure of minimal non-nilpotent groups
(also called \textit{Schmidt groups}) given by \cite{8}.

It is well-known that a finite nilpotent group $G$ contains elements of order $\exp⁡(G)$. Moreover, all sections of $G$ have this property. Under the above notation,
this can be written alternatively as
$$\varphi(S)\neq 0 \mbox{ for any section } S \mbox{ of } G.\0(1)$$Our main theorem shows that the converse is also true, that is we have the following nilpotency criterion.

\bigskip\noindent{\bf Theorem 1.} {\it Let $G$ be a finite group. Then $G$ is nilpotent if and only if $\varphi(S)\neq 0$ for any section $S$ of $G$.}
\bigskip

Note that (1) implies
$$\varphi(S)\neq 0 \mbox{ for any subgroup } S \mbox{ of } G\0(2)$$and in particular $$\varphi(G)\neq 0.\0(3)$$We observe that the condition (3) is not sufficient to guarantee the nilpotency of $G$, as shows the elementary example $G=\mathbb{Z}_6\times S_3$; we can even construct a non-solvable group $G$ for which $\varphi(G)\neq 0$, namely $G=\mathbb{Z}_n\times H$, where $H$ is a simple group of exponent $n$. A similar thing can be said about the condition (2).

\bigskip\noindent{\bf Example.} Let $G$ be a nontrivial semidirect product of a normal subgroup isomorphic to $$E(5^3)=\langle x,y \mid x^5=y^5=[x,y]^5=1, [x,y]\in Z(E(5^3))\rangle$$by
a subgroup $\langle a\rangle$ of order $3$ such that $a$ commutes with $[x,y]$. Then $G$ is a non-CLT group of order $375$, more precisely it does not have subgroups of order $75$. We infer that its subgroups are: $G$, all subgroups contained in the unique Sylow $5$-subgroup, all Sylow $3$-subgroups, and all cyclic subgroups of order $15$. Clearly, $G$ satisfies the condition (2), but it is not nilpotent.

Finally, we note that our criterion can be used to prove the non-nilpotency of a finite group by looking to its sections. In \cite{9} we have determined several classes of groups $G$ satisfying $\varphi(G)=0$, such as dihedral groups $D_{2n}$ with $n$ odd, non-abelian $P$-groups of order $p^{n-1}q$ {\rm(}$p>2,q$ primes, $q \mid p-1${\rm)}, symmetric groups $S_n$ with $n\geq3$, and alternating groups $A_n$ with $n\geq4$. These examples together with Theorem 1 lead to the following corollary.

\bigskip\noindent{\bf Corollary 2.} {\it If a finite group $G$ contains a section isomorphic to one of the above groups, then it is not nilpotent.}

\section{Proof of Theorem 1}

We will prove that a finite group all of whose sections $S$ satisfy $\varphi(S)\neq 0$ is nilpotent. Assume that $G$ is a counterexample of minimal order. Then $G$ is a Schmidt group since all its proper subgroups satisfy the hypothesis. By \cite{8} (see also \cite{2,7}) it follows that $G$ is a solvable group of order $p^mq^n$ (where $p$ and $q$ are different primes) with a unique Sylow $p$-subgroup $P$ and a cyclic Sylow $q$-subgroup $Q$, and hence $G$ is a semidirect product of $P$ by $Q$. Moreover, we have:
\begin{itemize}
\item[-] if $Q=\langle y\rangle$ then $y^q\in Z(G)$;
\item[-] $Z(G)=\Phi(G)=\Phi(P)\times\langle y^q\rangle$, $G'=P$, $P'=(G')'=\Phi(P)$;
\item[-] $|P/P'|=p^r$, where $r$ is the order of $p$ modulo $q$;
\item[-] if $P$ is abelian, then $P$ is an elementary abelian $p$-group of order $p^r$ and $P$ is a minimal normal subgroup of $G$;
\item[-] if $P$ is non-abelian, then $Z(P)=P'=\Phi(P)$ and $|P/Z(P)|=p^r$.
\end{itemize}We infer that $S=G/Z(G)$ is also a Schmidt group of order $p^rq$ which can be written as semidirect product of an elementary abelian $p$-group $P_1$ of order $p^r$ by a cyclic group $Q_1$ of order $q$ (note that $S_3$ and $A_4$ are examples of such groups). Clearly, we have $$\exp(S)=pq.$$On the other hand, it is easy to see that $$L(S)=L(P_1)\cup\{Q_1^x \mid x\in S\}\cup\{S\}.$$Thus, the section $S$ does not have cyclic subgroups of order $pq$ and consequently $\varphi(S)=0$, a contradiction. This completes the proof. \hfill\rule{1,5mm}{1,5mm}


\bigskip

\begin{thebibliography}{10}
\bibitem{1} {\bf H. Amiri, S.M. Jafarian Amiri, I.M. Isaacs}, {Sums of element orders in finite groups}, {\em Comm. Algebra} {\bf 37} (2009), no. 9, 2978--2980. 
\bibitem{2} {\bf A. Ballester-Bolinches, R. Esteban-Romero, D.J.S. Robinson}, {On finite minimal non-nilpotent groups}, {\em Proc. Amer. Math. Soc.} {\bf 133} (2015), 3455–-3462.
\bibitem{3} {\bf M. Garonzi, M. Patassini}, {Inequalities detecting structural properties of a finite group}, {\em Comm. Algebra} {\bf 45} (2017), no. 2, 677--687.
\bibitem{4} {\bf M. Herzog, P. Longobardi, M. Maj}, {An exact upper bound for sums of element orders in non-cyclic finite groups}, http://arxiv.org/abs/1610.03669.
\bibitem{5} {\bf T. De Medts, M. T\u arn\u auceanu}, {Finite groups determined by an inequality of the orders of their subgroups}, {\em Bull. Belg. Math. Soc. Simon Stevin} {\bf 15} (2008), 699–-704.
\bibitem{6} {\bf T. De Medts, T\u arn\u auceanu}, {An inequality detecting nilpotency of finite groups}, http://arxiv.org/abs/1207.1020.
\bibitem{7} {\bf V.S. Monakhov}, {The Schmidt subgroups, its existence, and some of their applications}, {\em Ukraini. Mat. Congr. 2001}, Kiev, 2002, Section 1, 81--90.
\bibitem{8} {\bf O.Yu. Schmidt}, {Groups whose all subgroups are special}, {\em Mat. Sb.} {\bf 31} (1924), 366--372.
\bibitem{9} {\bf M. T\u arn\u auceanu}, {A generalization of the Euler's totient function}, {\em Asian-Eur. J. Math.} {\bf 8} (2015), no. 4, article ID 1550087.
\end{thebibliography}
\end{document}